\documentclass[12pt]{article}
\usepackage{amsmath}
\usepackage{amssymb}
\usepackage{amscd}
\usepackage{amsfonts}
\newcounter{theorem}

\newcounter{theoremcounter}

\newcounter{remarkcounter}
\newcounter{propositioncounter}
\newcounter{definitioncounter}

\newcounter{corollarycounter}

\newtheorem{theorem}[theoremcounter]{Theorem}

\newtheorem{remark}[remarkcounter]{Remark}
\newtheorem{proposition}[propositioncounter]{Proposition}
\newtheorem{definition}[definitioncounter]{Definition}

\newtheorem{corollary}[corollarycounter]{Corollary}

\newcommand{\be}{\begin{equation}}
\newcommand{\ee}{\end{equation}}

\newcommand{\ber}{\begin{eqnarray}}
\newcommand{\eer}{\end{eqnarray}}

\newcommand{\refm [1]} {(\ref{#1})}
\newcommand{\nin}{\noindent}

\title{Nonstationary queues: Estimation of the rate of convergence }
\date{}
\author {{\bf Boris L. Granovsky}
\thanks{e-mail: mar18aa@techunix.technion.ac.il} \\
Department of Mathematics, Technion-Israel Institute of
Technology,\\ Haifa, 32000, Israel.\\
  {\bf A.
Zeifman}
\thanks{e-mail:zai@uni-vologda.ac.ru} \\
Vologda State Pedagogical University and\\ Vologda Science
Coordination Centre CEMI RAS. Vologda, Russia.}
\begin{document}
\maketitle
\vskip 5cm
\newpage

\section{Introduction.}

\nin It is commonly acknowledged that explicit expressions for the
transient behaviour of functionals of stochastic models can be
found only in a few special cases. In view of this, efforts of
generations of probabilists have been focused on the study of the
 rate of convergence, as time
$t\to \infty ,$ to the steady state of a process. During the last
two decades a remarkable progress was made in this field,
regarding time-homogeneous Markov chains, as a result of
implementation and development of sophisticated techniques:
coupling, logarithmic Sobolev inequalities, the Poincar\'e
inequality and its versions arising from the variational
interpretation of eigenvalues, and duality. For reference, we
recommend recent review and research papers \cite{ald}, \cite{ch},
\cite{doorn85}, \cite{fill}, \cite{Diac1}, \cite{doo00},
\cite{doo01}. The new stream of nowadays research was motivated by
new fields of applications, s.t. algorithms of Monte Carlo for
simulation of Markov chains and enumeration algorithms in computer
science and group theory.

\nin  A variety  of problems for queueing models, among them
stability, were investigated by V. Kalashnikov \cite{vk}.

\noindent
 There is also a growing interest in time-nonhomogeneous Markov chains
(see  \cite{gio87}, \cite{hey84}, \cite{kel82}) that  model a
large number of real queuing systems. Most of research on
nonhomogeneous queues is devoted to various methods of
approximation of their transient behaviour(see \cite{dicres95},
\cite{man95}, \cite{mas94} and references therein.) A quite
different approach is used in the study of annealing processes
(for references see \cite{dm}). Here Markov chains considered have
exponentially vanishing intensities, as $t\to \infty.$ The
research in this field is aimed to estimate the rate of
convergence of the process to the invariant distribution.

\noindent Our paper is devoted to the estimation of the rate of
 of exponential convergence of
 nonhomogeneous queues exhibiting different types of ergodicity.
  The main tool of our study is the method,
which was proposed by the second author in the late 1980s (see
\cite{z85}) and was subsequently extended and developed in
different directions in a series of joint papers by the authors of
the present paper. The method originated from the idea of Gnedenko
and Makarov(\cite{gm}) to employ the logarithmic norm of a matrix
to the study of the problem of stability of Kolmogorov system of
differential equations associated with nonhomogeneous Markov
chains. The method is based on the following two ingredients: the
logarithmic norm of a linear operator and a special similarity
transformation of the  matrix of intensities of the  Markov chain
considered.

\nin In the paper we apply the method to a new general class of
 Markov queues with a special form of nonhomogenuity that is common
in applications. Namely, we consider the case of asymptotically
periodic arrival and service rates.

\nin  The present paper is a substantially  modified and extended
version of our
 preliminary report \cite{gz00}.

\section{Preliminaries.}
\nin  We  consider a contionuous time nonhomogeneous birth and
 death process (BDP)

\nin  $X\left(t\right),\ t\ge 0 $
  on the state space $E=\left\{ 0,1,\ldots ,N\right\}, $ $N\leq \infty, $
  with time dependent intensities of birth
   $\lambda _n\left( t\right),\ t\ge 0, $
   and death $\mu _n\left( t\right), \ t\ge 0,\ n\in E. $
\nin Namely,

\begin{equation}
P(X(t+h) =j\vert X(t)=i)= \left\{
\begin{array}{cc}
\lambda _i\left( t\right) h+o\left( h\right), & {if } \quad j=i+1
\\ \mu _i\left( t\right)  h+o\left( h\right), & \mbox {if }
\quad j=i-1
\\
1-\left( \lambda _i\left( t\right) +\mu _i\left( t\right) \right)
 h+o\left( h\right), & \mbox {if } \quad j=i \\ o\left( h\right),
& \mbox {if } \quad \left| i-j\right| >1,
\end{array}
\right.  \label{eq111}
\end{equation}

\nin where $h>0$ and, for the sake of brevity, $o(h)=o_i(t,h), \
i\in E$ denotes different quantities s.t. for all $t\ge 0$, $$
\lim_{h\to 0} \frac{\sup_i |o_i(t,h)|}{h}=0. $$

\nin  Denote

\begin{equation} p_{ij}\left( s,t\right) =Pr \left( X\left(
t\right) =j\mid X\left( s\right) =i\right), \ i,j\in E, \quad
0\leq s\leq t
\end{equation}
and $p_i\left( t\right) =\Pr \left( X\left( t\right) =i\right), \
i\in E, \quad t\ge 0 $ the transition and the state probabilities
respectively.

\nin Let ${\bf p}={\bf p}\left( t\right) = (p_0(t),\ldots,
p_N(t))^T, \ t\ge 0$
 be a column vector of probabilities of states and let $A(t), \quad t\ge
0 $ be the intensity matrix induced by
$\left(\ref{eq111}\right)$.

\noindent Then the evolution of the process  $X(t)$ is described
by the forward Kolmogorov system

\begin{equation}
\frac{d{\bf p}}{dt}=A\left( t\right) {\bf p},\quad {\bf p}={\bf p}(t),
\quad t\ge 0.
\label{eq112}
\end{equation}

\nin The Cauchy operator $U\left( t,s\right),\quad 0\leq s\leq t $
of  (\ref{eq112}) is given by the matrix $U^T\left( t,s\right)
=P\left(s,t\right):= \left( p_{ij}\left( s,t\right) \right)
_{i,j=0}^N.$  We denote  throughout the paper by $\Vert \bullet
\Vert$ the $l_1-$ norm, i. e.  $\left\| {\bf x} \right\|
=\sum_{i\in E} \left| x_i\right|,$  for ${\bf x=}\left( x_0,\ldots
,x_N\right)^T$ and $\left\| B\right\| = \sup_{j\in
E}\sum\limits_{i\in E}\left| b_{ij}\right| ,$ for $ B=\left(
b_{ij}\right) _{i,j=0}^N.$

\noindent Let $\Omega $ be the set of vectors $\Omega $= $\{{\bf
x=}\left( x_0,\ldots ,x_N\right) ^T:{\bf x}\geq 0,\quad \left\|
{\bf x}\right\| =1\}.$  Assuming  that $ {\bf p}(t)\in \Omega, \
t\ge 0$ is a solution of (\ref{eq112}), we consider the following
types of ergodicity of continuous time Markov chains.

\begin{definition}.
A Markov chain $X(t),\quad t\ge 0$ on $E$ is called

\noindent {\bf ergodic}, if there exists a probability measure
$\pi=(\pi(0),\ldots,\pi(N))^T$   on $E,$

\nin such that
$\lim\limits_{t\rightarrow \infty }\left\| {\bf p}\left( t\right) -%
{\bf \pi }\right\| =0,$  for all $ {\bf p}\in \Omega;$

\noindent {\bf weakly ergodic}, if $\lim\limits_{t\rightarrow
\infty }\left\| {\bf p}^{(1)}\left( t\right) -{\bf p}^{(2)}\left(
t\right) \right\| =0,$ for all ${\bf p}^{(i)}\in \Omega, \ i=1,2;$
\nin

\nin {\bf null-ergodic }, if $p_n(t) \to 0,\quad n\in E,$ as
 $t\rightarrow \infty ,$ for all
 ${\bf p}\in \Omega$ and

\nin {\bf quasi-ergodic}, if there exists a probability measure
 ${\bf \tilde \pi }$ on $E,$ such that
$$\lim_{t \to \infty } \left\| \frac{1}{t} \int_{0}^{t} {\bf
p}\left( \tau\right)\,d\tau - {\bf \tilde \pi }\right\| = 0,$$ for
all ${\bf p}\in \Omega  .$

\end{definition}

\begin{remark}.
BDP can be both weakly ergodic and null-ergodic,
 see Example 1 in \cite{gz00}. Hovewer, quasi-ergodic  BDP can not
  be null-ergodic.
\end{remark}

\noindent We deal with the class of nonhomogeneous BDP given by

\begin{equation}
\lambda _n\left( t\right) =\lambda _n a\left( t\right),\
\mu_n\left( t\right) =\mu _n
 b\left( t\right) ,\quad t\ge 0 \quad
n\in E.
 \label{eq101}
\end{equation}
\nin We assume throught the  paper  that the basic functions
$a\left( t\right)\ge 0,\ b\left( t\right)\ge 0,\ t\ge 0$ are
locally integrable on $\left[ 0,\infty \right) $ and
asymptotically periodic with periods $T_1, T_2$ correspondingly.
The latter means that

\begin{equation}
a ( t) = a_{1}(t)+ a_{2}(t), \ b ( t) = b_{1}(t)+ b_{2}(t),\quad
t\ge 0, \label{ap1}
\end{equation}

\nin where
\begin{equation}
a_{1}(t+T_{1}) = a_{1}(t), \ b_{1}(t +T_{2}) = b_{1}(t) , t \ge 0,
\quad \lim_{t \to \infty} \Big(|a_{2}(t)| + |b_{2}(t)|\Big) = 0.
\label{ap2}
\end{equation}

\nin Regarding the rates $\lambda_n, \mu_n,$ we assume that  \be
\lambda _n>0,\quad n=0,\dots,N-1,\quad \mu _n>0,\ n=1,\ldots,N,\
\mu_0=0. \end{equation} We also assume that in the case $N=\infty$
there exist the limits \be \lim\limits_{n\rightarrow \infty
}\lambda _n=\lambda >0,\ \lim\limits_{n\rightarrow \infty }\mu
_n=\mu>0, \label{eq102}
\end{equation}
\nin while in the case  $N<\infty$,

\begin{equation}
 \lambda_N=0,\quad
\lambda: =\lambda _{N-1},\ \mu: =\mu _N.  \label{eq103}
\end{equation}

\noindent Some particular cases of the above  setting were studied
in \cite{gz1}, \cite{ki92}, \cite{g}. In the aforementioned papers
only ergodic BDP's were treated. In the present paper  we shall
investigate the rate of convergence in other cases given by
Definition 2.

\nin For reader's convenience, we recall the method of bounding of
the
 the rate of convergence of BDP's that was introduced  in \cite{z95b}.

\noindent First, we   recall the definition of the logarithmic
norm that was proposed for finite-dimensional spaces  by Lozinskij
\cite{L} and generalized to  Banach spaces by Daleckij and Krein
\cite{DK}.

\begin{definition}
Let $B\left( t\right) ,\ t\ge 0$ be a one-parameter family of
bounded linear operators on a Banach space ${\cal B}$ and let $I$
denote the identity operator.

\noindent For a given $t\ge 0$, the number

\begin{equation}
\gamma \left( B\left( t\right) \right) =\lim\limits_{h\rightarrow
+0}\frac{%
\left\| I+hB\left( t\right) \right\| -1}h  \label{eq116}
\end{equation}
is called the logarithmic norm of the operator $B\left( t\right) .$
\end{definition}

\noindent If ${\cal B}$ is a $(N+1)-$dimensional vector space with
$l_1$- norm, so that the operator $B(t)$ is given by the matrix
$B(t)=\left( b_{ij}(t)\right) _{i,j=0}^N,\quad t\ge 0$, then the
logarithmic norm of $B(t)$ can be found explicitly:

\begin{equation}
\gamma \left( B\left( t\right) \right) =\sup\limits_j\left( b_{jj}\left(
t\right) +\sum\limits_{i\neq j}\left| b_{ij}\left( t\right) \right| \right)
,\quad t\ge 0.  \label{eq1000}
\end{equation}

\noindent Associate now the family of operators $B(t), \quad t\ge 0$ with
the system of differential equations

\begin{equation}
\frac{d{\bf x}}{dt}=B\left( t\right) {\bf x},\quad t\ge 0,
\label{eq1111}
\end{equation}
where the functions $b_{ij}(t),\ 0\le i,j\le N$ are assumed to be
locally integrable on
$%
[0,\infty).$ Then the  logarithmic norm of the operator $B(t)$ is
related to  the Cauchy operator $V(t,s), \ 0\le s\le t$ of the
system $\left(\ref {eq1111}\right):$

\begin{equation}
\gamma \left(B\left( t\right) \right) =\lim\limits_{h\rightarrow
+0}\frac{%
\left\| V\left( t+h,t\right) \right\| -1}h, \quad t\ge 0.
\label{eq1100}
\end{equation}

\noindent From the latter  one can  deduce the  following bounds
on the ${\cal B} $-norm of the Cauchy operator $V(t,s),\quad 0\le
s\le t$:

\begin{equation}
e^{-\int\limits_s^t\gamma \left( -B\left( \tau \right) \right) \ d\tau
}\leq
\left\| V\left( t,s\right) \right\| \leq e^{\int\limits_s^t\gamma \left(
B\left( \tau \right) \right) \ d\tau },\quad 0\le s\le t.  \label{eq1117}
\end{equation}

\nin Moreover, for any solution ${\bf x}(t)\in {\cal B},\quad t\ge
0$ of $\left(\ref {eq1111}\right)$ we have

\begin{equation}
\left\| {\bf x}\left( t\right) \right\| \geq
e^{-\int\limits_s^t\gamma \left( -B\left( \tau \right) \right) \
d\tau } \left\| {\bf x}\left( s\right) \right\|. \quad
\label{eq1118}
\end{equation}

\noindent We will also make use of the fact that if ${\cal B}=l_1$
and all non-diagonal elements of $B$
\begin{equation}
b_{ij}\left( t\right) \geq 0, \quad i\neq j,\quad t\geq 0,
\label{eq1114}
\end{equation}
then, by $\left(\ref{eq1000}\right),$
\begin{equation}
\gamma \left( B\left( t\right) \right)
=\sup\limits_j\sum\limits_ib_{ij}\left( t\right), \quad t\ge 0,
\end{equation}

\noindent and, consequently, for any solution ${\bf x}\left(
t\right), \ t\ge 0$ of $\left(\ref{eq1111}\right),$ s.t. ${\bf
x}\left( s\right)\geq {\bf 0},$ we have

\begin{equation}
\left\| {\bf x}\left( t\right) \right\| \geq
e^{\int\limits_s^t\inf\limits_j\sum\limits_ib_{ij}\left( \tau
\right) \ d\tau } \left\| {\bf x}\left( s\right) \right\|, \quad
0\le s \le t . \label{eq1119}
\end{equation}

\nin  We adopt further on  the  definition of the decay function
$\beta(t)$ of a nonhomogeneous ergodic Markov chain,  given by the
authors in \cite{gz1}. (Note  that in the case of a homogeneous
Markov chain, $\beta(t)\equiv\beta,\ t\ge 0$ is the spectral gap.)

\noindent As a result of the previous discussion,  we derive
 the following two-side bounding of the
decay function $\beta (t),\ t\ge 0$, provided that  the matrix
$B(t), \ t\ge 0$ obeys  the condition
$%
\left( \ref{eq1114}\right): $
\begin{equation} \underline
c(t):=\inf\limits_j\sum\limits_ib_{ij}\left( t\right) \le -
\beta(t)\le \bar c(t)
:=\sup\limits_j\sum\limits_ib_{ij}\left( t\right) ,\quad
t\ge 0.  \label{eq2114}
\end{equation}

\noindent In the  method considered, the inequalities
(\ref{eq1117}), (\ref{eq1118}), (\ref{eq2114}) constitute  the
main tool for bounding the decay function, as well as the rates of
convergence of  other types of Markov chains in Definition 1.
Finally, we observe that in the case $%
\underline c(t)=\bar c(t),\ t\ge 0,$ the estimate $\left(\ref{eq2114}%
\right) $ is exact.

\noindent A straightforward implementation of $\left(\ref{eq2114}%
\right)$ to the operator $A(t)$ is not effective, since, by the
definition of the intensity matrix and
$\left(\ref{eq1000}\right)$, we have $\gamma \left(A\left(
t\right) \right)=0,\ t\ge 0.$

\noindent For this reason, we rewrite $\left( \ref{eq112}\right) $
as a  system of nonhomogeneous differential equations with a
nonsingular matrix. Namely, substituting in $\left(
\ref{eq112}\right) $

\begin{equation}
p_0\left( t\right) =1-\sum\limits_{i\ge 1} p_i\left( t\right), \quad t\ge
0,
\end{equation}
we have

\begin{equation}
\frac{d{\bf z}}{dt}= B(t){\bf z}(t) + {\bf f}(t) ,\quad t\ge 0,
\label{eq304}
\end{equation}
where
\begin{eqnarray}
 B(t)= \{a_{ij}(t) - a_{i0}(t), \quad i,j=1,\ldots, N\},
\nonumber
\\
{\bf z}(t)=( p_1(t),\ldots , p_N(t)) ^T,\quad  \nonumber \\
{\bf f}(t) =\left( a_{10}(t) ,\ldots ,a_{N0}(t) \right) ^T,\quad t\ge 0.
\label{eq8}
\end{eqnarray}

\noindent The solution of $\left( \ref{eq304}\right) $ is given by

\begin{equation}
{\bf z}\left( t\right) =V\left( t,s\right)  {\bf z}\left( s\right)
+\int\limits_s^tV\left( t,\tau \right)  {\bf f}\left( \tau \right)
\,d\tau, \quad 0\le s \le t,  \label{eqS5}
\end{equation}
where $V\left( t,s\right) $ is the Cauchy operator of
the system $\left( \ref {eq304}\right) .$

\noindent  The following simple relationship holds between pairs
of solutions

\[
{\bf z}^{(i)}={\bf z}^{(i)}(t), \quad {\bf p}^{(i)}= {\bf p}%
^{(i)}(t),\quad t\ge 0, \quad i=1,2
\]
of $\left(\ref{eq304}\right) $ and $\left(\ref{eq112}\right)$ respectively:

\begin{eqnarray*}
\left\| {\bf p}^{(1)}-{\bf p}^{(2)}\right\| &=&\left|
p_0^{(1)}-p_0^{(2)} \right| +\sum\limits_{i\geq 1}\left|
p_i^{(1)}-p_i^{(2)}\right| =\left| 1-\sum\limits_{i\geq 1}p_i^{(1)}-\left(
1-\sum\limits_{i\geq 1}p_i^{(2)}\right) \right| \\
+\left\| {\bf z}^{(1)}-{\bf z}^{(2)}\right\| &=&\left|
\sum\limits_{i\geq 1}\left( p_i^{(1)}-p_i^{(2)}\right) \right| +\left\|
{\bf z}^{(1)}-{\bf z}^{(2)}\right\|, \quad t\ge 0.
\end{eqnarray*}

\noindent Consequently,
\begin{equation}
\left\| {\bf z}^{(1)}-{\bf z}^{(2)}\right\| \leq \left\| {\bf p}%
^{(1)} -{\bf p}^{(2)}\right\| \leq 2\left\| {\bf z}^{(1)}-{\bf z}%
^{(2)}\right\|, \quad t\ge 0,  \label{eqS4}
\end{equation}
\noindent which will be used in our subsequent  study.

\noindent At this point, we  note that the matrix $ B(t),\quad
t\ge 0$ in $\left( \ref {eq8}\right) $ lacks the property $\left(
\ref{eq1114}\right)$ of the original intensity matrix $A(t), \
t\ge 0 $ in $\left( \ref{eq112}\right) $. This fact  prevents the
implementation of $\left( \ref{eq2114}\right) .$ \noindent To
overcome this  obstacle, we employ a similarity transformation of
$B(t),\ t\ge 0,$ that restores the  property $\left(
\ref{eq1114}\right)$. In \cite{gz2} it was proven the existence of
such transformation for general finite homogeneous Markov chains.
However, its explicit construction is known so far for homogeneous
BDP's only. The required transformation which is in the core of
the method considered, was  suggested in \cite{z89}-\cite{z91b}.
It is given by the upper triangular $N\times N$ matrix

\begin{equation}
D=\left(
\begin{array}{cccc}
d_0 & d_0 & d_0 & \ldots \\
0 & d_1 & d_1 & \ldots \\
0 & 0 & d_2 & \ldots \\
\ldots & \ldots & \ldots & \ldots
\end{array}
\right)  \label{eqS7}
\end{equation}
where $d_i>0, \quad i=0,1, \ldots, N-1. $ We have
\[
D^{-1}=\left(
\begin{array}{ccccc}
d_0^{-1} & -d_1^{-1} & 0 & \ddots &  \\
0 & d_1^{-1} & -d_2^{-1} & 0 & \ddots \\
\ddots & 0 & \ddots d_2^{-1} & \ddots & \ddots \\
& \ddots & \ddots & \ddots & \ddots \\
&  & \ddots & \ddots & \ddots
\end{array}
\right) .
\]

\noindent  Applying this transformation to the matrix $B(t), \ t\ge 0$
 in $\left(%
\ref{eq8}\right),$  leads to the matrix $DB(t)D^{-1}, \ t\ge 0$
with the desired property $\left( \ref{eq1114}\right): $

\begin{equation}
DB(t)D^{-1}=\left(
\begin{array}{ccccc}
-\left( \lambda _0(t)+\mu _1(t)\right) & d_0 d_1^{-1} \mu _1(t) & 0 &
& \ldots \\ d_1 d_0^{-1} \lambda _1(t) & -\left( \lambda _1(t)+\mu
_2(t)\right) & d_1 d_2^{-1} \mu _2(t) & 0 &
\\ 0 & d_2 d_1^{-1} \lambda _2(t) & \ddots & \ddots &
\ddots \\ \vdots & 0 & \ddots & \ddots & \ddots \\ & \ldots &
\ddots & \ddots & \ddots
\end{array}
\right) . \label{eqS8}
\end{equation}

\noindent  We set further on, $d_0=1,$ and $d_N=0$ for $N < \infty
,$ and denote $\delta_{i+1}=d_{i+1}/d_i>0,\quad i= 0,\ldots ,
N-1.$ By $\left(\ref{eq1000}\right),$ the logarithmic norm of the
above matrix equals to

\begin{equation}
\gamma (DB(t)D^{-1})=-\inf_k\alpha _k\left( t\right) ,\quad t\ge 0,
\label{eq21160}\end{equation}

\nin where

\begin{equation}
\alpha _k\left( t\right) =\lambda _k\left( t\right) +\mu
_{k+1}\left( t\right) -\delta _{k+1} \lambda _{k+1}\left( t\right)
-\delta _k^{-1} \mu _k\left( t\right) , \quad t\ge 0, \quad
k=0,\ldots ,N-1.  \label{eq2116}
\end{equation}
 \nin (Here and in what follows $\delta_0=0$.)

\noindent Next, the substitution ${\bf z}(t)=D^{-1}{\bf y}(t),\
t\ge 0,$ transforms $\left( \ref{eq304}\right) $ into the system
with the matrix $DB(t)D^{-1},\ t\ge 0.$ Following \cite{z89},
\cite{z91a}, we denote
\begin{equation}
\left\| {\bf z}\right\| _{1D}=\left\| D{\bf z}\right\| ,\quad t\ge 0,
\label{eq2120}
\end{equation}
and $g=\inf_{k>0}d_k,$ to obtain, by virtue of $\left(
\ref{eqS7}\right) ,$

\begin{equation}
\left\| {\bf z}\right\| _{1D}\ge \frac{g}{2}\left\| {\bf z}\right\|,
\quad t\ge 0.  \label{eq2117}
\end{equation}

\noindent Now we apply $\left(\ref {eq2114}\right)$
to  the matrix $DB(t)D^{-1}, \ t\ge 0,$ to derive 
 the following effective way of two-side bounding of the decay
 function $\beta(t), \ t\ge 0$ of  non-homogeneous BDP's:
\begin{equation}
\underline \alpha (t) \le \beta (t) \le \bar \alpha (t), \quad t \ge 0,
\label{eq4118}
\end{equation}
where
\begin{equation}
\underline \alpha(t)=\inf_k \alpha_k(t), \quad \bar
\alpha(t)=\sup_k\alpha_k(t), \quad t\ge 0 \label{eq4119}
\end{equation}

\nin In this connection note that, by the definition in
\cite{gz1}, the decay function $\beta(t),\ t\ge 0$ should satisfy
$\int\limits_0^{+\infty} \beta \left( t\right)dt=+\infty,$ (but is
not required to be nonnegative for all $t\ge 0$).

\nin  Therefore, $\left( \ref{eq4118} \right)$ make sense provided
the weights $\delta_k,\ k=0,\ldots,N-1$ in $\left(
\ref{eq2116}\right) $ are s.t. $\int\limits_0^{+\infty} \underline
\alpha(t)dt=+\infty.$  The way of finding such weights is
discussed further on in the present paper.

\noindent Next, with the help of $\left(\ref{eq1117}\right)$,
$\left(\ref{eq2117}\right)$ and  $\left(\ref{eqS4}\right)$ we
obtain for weakly ergodic BDP's
:

\begin{eqnarray}
\left\| {\bf p}^{(1)}\left( t\right) -{\bf p}^{(2)}\left( t\right)
\right\|_{1D} \le e^{-\int\limits_s^t\underline \alpha \left(
u\right) du} \left\| {\bf p}^{(1)}\left( s\right) -{\bf
p}^{(2)}\left( s\right) \right\| _{1D}, \quad 0\le s\le t.
\label{eq21180}
\end{eqnarray}

\nin and

\begin{eqnarray}
\left\| {\bf p}^{(1)}\left( t\right) -{\bf p}^{(2)}\left( t\right)
\right\| \le \frac 4g e^{-\int\limits_s^t\underline \alpha \left(
u\right) du} \left\| {\bf z}^{(1)}\left( s\right) -{\bf
z}^{(2)}\left( s\right) \right\| _{1D}, \quad 0\le s\le t.
\label{eq2118}
\end{eqnarray}
\noindent Here, in accordance with our notation
$\left(\ref{eq2120}\right),$
\begin{equation}
\left\| {\bf z}^{(1)}\left( s\right) -{\bf z}^{(2)}\left( s\right)
\right\| _{1D} \le \sum\limits_{i\geq 1}q_i \left| p_i^{(1)}\left(
s\right) -p_i^{(2)}\left( s\right) \right|, \quad t\ge 0,
\label{eq2121}
\end{equation}
\noindent where $q_i=\sum\limits_{m=0}^{i-1}d_m,\quad i=1,\ldots,N.$

\noindent Clearly, (\ref{eq2118}), (\ref{eq2121})  are valid  also
for ergodic BDP's, by taking ${\bf p}^{(2)}\left( t\right) ={\bf
\pi },\quad t\ge 0.$

\nin The reasoning analogous to the preceding  one
 enables also to derive  a lower bound for
 $\left\| {\bf p}^{(1)}\left( t\right) -{\bf p}^{(2)}\left( t\right)
\right\|.$ For this purpose,
 for  a sequence of positive weights
$\mbox{\boldmath $\delta$}=
( \delta _1, \ldots, \delta_{N-1}),$
define
the quantities

\begin{equation}
\zeta _k\left( t\right) =\lambda _k\left( t\right) +\mu
_{k+1}\left( t\right) +\delta _{k+1} \lambda _{k+1}\left( t\right)
+\delta _k^{-1} \mu _k\left( t\right), \quad t\ge 0, \quad
k=0,\dots , N-1.  \label{eqf2}
\end{equation}

\nin Denoting
\begin{equation}
\bar\zeta \left( t\right) =\sup\limits_k\zeta _k\left( t\right),
\quad t\ge 0 .  \label{eqf3}
\end{equation}

\nin we then have, $\gamma (-DB(t)D^{-1}) =-\sup_k\zeta _k\left(
t\right) ,\quad t\ge 0,$  $(\left(\ref{eq1000}\right).$
Consequently, with the help of $\left(\ref{eqS4}\right),$ we
obtain

\begin{eqnarray}
&&\left\| {\bf p}^{(1)}\left( t\right) - {\bf p}^{(2)}\left(
t\right) \right\| \geq \left\| {\bf z}^{(1)}\left( t\right) -{\bf
z}^{(2)}\left( t\right) \right\| \geq  \nonumber   \\ &&\left\|
D\right\| ^{-1} \left\| D^{-1}\right\| ^{-1}
e^{-\int\limits_s^t\bar\zeta \left( u\right) du} \left\| {\bf
z}^{(1)} \left( s\right) -{\bf z}^{(2)}\left( s\right) \right\|
,\quad 0\le s\le t. \nonumber \\ \label{eqf10}
\end{eqnarray}
\vskip .5cm

\nin Now we are in a position to describe the organization of the
paper and its main results. As we already mentioned, we focus our
study on a class of nonhomogeneous BDP's with intensities (
\ref{eq101})-(\ref{ap2}) .
 Our main effort is
 devoted to estimation of $\underline{\alpha }(t),\ t\ge 0$
for these BDP's. This is done in Section 3 and 4, by analyzing
 the specific form of the expressions  $\left( \ref{eq2116}\right) $ in the case
considered. As a result, we derive estimates of the rate of
exponential convergence in different types of  ergodicity.
 As a by-product of this study, we obtain also estimates of
 important for applications
characteristics of transient behaviour of the above processes.
 The rest four sections are devoted to four particular
 nonhomogeneous queues widely known in the literature.
 We achieve an improvement of known estimates of the rate
 of convergence. We also obtain estimates of the expected
 length of queues. In the conclusion note that the stability properties of BDP's
  were explored
  in \cite{z85}, \cite{z94}, \cite{gz00} and \cite{z98}.

\section{Quasi-ergodic BDP}

 \nin Denote
\begin{equation}
f_k=\frac{\Delta\lambda \mu _{k+1}-\lambda _{k-1}\mu +\sqrt{\left(
\Delta\lambda \mu _{k+1}-\lambda _{k-1}\mu \right) ^2+4\Delta\lambda \mu
  \left( \lambda _{k-1}\mu _{k+1}-\lambda _k\mu _k\right) }}{%
2\Delta\lambda \mu },  \quad k=0,\ldots, N-1, \label{eq203}
\end{equation}
where $\Delta>0$ is a given number, and $\lambda, \mu$ are defined
in $\left( \ref{eq102}\right) $, $\left( \ref{eq103}\right). $

We denote $\tilde E=\left\{ 0,1,\ldots ,N-1\right\}, $ if
$N<\infty $ and $\tilde E=E,$ if $N=\infty $ and set in
\refm[ap1],\refm[ap2]

\begin{equation}
a_{m}: = \lim_{t \to \infty} \frac{1}{t} \int\limits_0^t a(u)\, du
= \frac{1}{T_{1}} \int\limits_0^{T_{1}} a_{1}(u)\, du; \quad
b_{m}: = \lim_{t \to \infty} \frac{1}{t} \int\limits_0^t b(u)\, du
= \frac{1}{T_{2}} \int\limits_0^{T_{2}}   b_{1}(u)\, du.
\label{eq2040}
\end{equation}

\begin{theorem}

\nin  Let for some $\Delta>1$ the following conditions fulfilled:

\noindent a) for any $k \quad \lambda_{k-1}\mu_{k+1}-\lambda_k\mu_k \ge 0;$

\noindent b)
\begin{equation}
\mu b_{m} - \Delta  \lambda a_{m} > 0, \label{eq2052}
\end{equation}

\noindent c)

\begin{equation}
\inf_{k\in \tilde E}f_k: = f > 0.  \label{eq205}
\end{equation}

Let $c \in (0;1)$ be arbitrary number such that a) $c < \frac{\mu_{k+1}}{\mu}, k \ge 0$ and b) $c \le f$. Then there exists a sequence $\mbox{\boldmath $\delta$}=(\delta_k,\quad
k\in \tilde E)>{\bf 0},$
 such that $\lim \inf \delta _k>1, \quad {if} \quad N=\infty $
and

(i)
\begin{equation}
\underline \alpha \left( t\right) \ge l(t): = c  \left( \mu
b\left( t\right) -\Delta\lambda a\left( t\right) \right),  \quad
t\ge 0 ,   \label{eq206}
\end{equation}

\nin (ii) BDP $X(t),\quad t\ge 0$ is weakly ergodic, so that
\begin{equation}
\left\| {\bf p}^{(1)}\left( t\right) -{\bf p}^{(2)}\left( t\right) \right\|_{1D}
\le e^{-\int\limits_s^t l \left( u\right) du}
\left\| {\bf p}^{(1)}\left( s\right) -{\bf p}^{(2)}\left( s\right) \right\|_{1D}
, \quad 0\leq s\leq t,   \label{eqS1001}
\end{equation}

\nin and

\begin{equation}
\left\| {\bf p}^{(1)}\left( t\right) -{\bf p}^{(2)}\left( t\right) \right\|
\leq \frac 4g  e^{-\int\limits_s^t l \left( u\right) du}
\sum\limits_{i\geq 1}q_i\left| p_i^{(1)}\left( s\right) -p_i^{(2)}
\left( s\right)
\right|, \quad 0\leq s\leq t,   \label{eqS1002}
\end{equation}

\nin where $g$ and $q_i$ are defined in $\left( \ref{eq2117}\right),$
$\left( \ref{eq2118}\right) $,

(iii) BDP $X(t),\quad t\ge 0$ is quasi-ergodic.

\end{theorem}

\vskip .5cm

\nin {\bf Proof.} The condition $0 < c \le f \le f_k$ implies the inequality

\begin{equation}
c^2\Delta\lambda \mu -c  \left( \Delta\lambda \mu _{k+1}-\lambda
_{k-1}\mu \right) +\lambda _k\mu _k-\lambda _{k-1}\mu _{k+1}\leq 0,
\quad k\in \tilde E
\label{eq210}
\end{equation}
\noindent for all $k$. Therefore using the inequality $c<  \frac{\mu_{k+1}}{\mu}$, we obtain

\begin{equation}
\frac{\mu_k}{\mu _{k+1}-c\mu }\leq \frac{\lambda_{k-1}+c\Delta\lambda }{%
\lambda_k}, \quad k\in \tilde E.  \label{eq209}
\end{equation}

Put now

\begin{equation}\delta _k \in \left[\frac{\mu_k}{\mu_{k+1}-c\mu }, \frac{\lambda_{k-1}+c\Delta\lambda }{\lambda_k} \right], \quad k\in \tilde E.
\label{eq207}
\end{equation}

Then

\begin{equation}
\delta_{k+1}  \lambda_{k+1}-\lambda_k\leq c\Delta\lambda, \quad
\mu_{k+1}-\delta_k^{-1}  \mu_k\geq c\mu, \quad k\in \tilde E,
\label{eq206,5}
\end{equation}

\noindent and, finally

\begin{equation}
\alpha_k\left( t\right) \geq c  \left( \mu b\left( t\right)
-\Delta\lambda a\left( t\right) \right),
\quad k\in \tilde E,\quad t\ge 0,
\label{eq1012}
\end{equation}

\nin because in the case considered

\begin{equation}
\alpha_k \left(t\right) = b \left( t\right)   \left( \mu_{k+1}-
\delta_k^{-1}  \mu_k\right) -a\left( t\right)
\left( \delta_{k+1}
\lambda_{k+1}-\lambda_k\right) , \quad k\in \tilde E,\quad t\ge 0.
\label{eq1011}
\end{equation}

Now, for the case  $N=\infty$ we get from (\ref{eq207}) as $k \to \infty$

\[
\liminf \delta_k\geq \frac \mu {\mu -c\mu }=\frac 1{1-c}>1
\]
\noindent because $c<1$.

\vskip .5cm

\begin{corollary}.
If under the assumptions of Theorem 1, the  BDP $X(t),\quad t\ge
0$ considered  has a stationary distribution ${\bf \pi }$ on $E, $
then the BDP is ergodic and

\begin{equation}
\left\| {\bf p}\left( t\right) -{\bf \pi }\right\| \leq
\frac 4g  e^{-\int\limits_s^t l \left( u\right) du}
\sum\limits_{i\geq 1}q_i\left|  p_i\left( s\right) -\pi_i
\right|, \quad 0\le s\le t .
\label{equS12}
\end{equation}
\end{corollary}
\vskip .5cm

\vskip .5cm

\begin{corollary}.
Under the assumptions of Theorem 1, the  BDP $X(t),\quad t\ge 0$
considered is exponentially weakly ergodic and for any
$\varepsilon
> 0$ there exists a constant $K$ such that

\begin{equation}
exp\left(-\int\limits_s^t l \left( u\right) du \right) \le   K e^{-(l-\varepsilon)(t-s)},
\quad t\ge 0,
\label{eq2101}
\end{equation}

\nin where $l = c\left(\mu b_{m} - \Delta  \lambda a_{m}\right).$

\end{corollary}
\vskip .5cm

\begin{theorem}
Under the assumptions of Theorem 1, for any $\varepsilon > 0$
there exists a constant $K$ such that
\begin{equation}
\Pr \left( X\left( t\right) \leq j\mid X\left( 0\right) =0\right) \geq
1 - q_{j+1} ^{-1} K\lambda_{0} a_{m} \left( 1+\frac{e^{(l-\varepsilon)T_{1}}}
{e^{(l-\varepsilon)T_{1}} - 1} \right).
\label{eqS16}
\end{equation}
\end{theorem}

\vskip .5cm

\nin {\bf Proof.} Put in $\left( \ref{eq304}\right) $ ${\bf y}\left( t\right) =D{\bf z}\left( t\right) $, then we have

\begin{equation}
\frac{d{\bf y}}{dt}=DB\left( t\right) D^{-1}{\bf y}+D{\bf f}\left(
t\right),\quad t\ge 0.
\label{eqS14}
\end{equation}

The matrix $DB\left( t\right) D^{-1}$ satisfies   $\left(
\ref{eq1114}\right)$, so under the  conditions ${\bf y}\left(
0\right)= {\bf 0}, \quad {\bf f}(t) \geq {\bf 0}, \quad t\ge 0,$
we get ${\bf y}\left( t\right) \geq {\bf 0}, \quad t\geq 0.$ Let
$V_1\left( t,s\right) $ be a Cauchy matrix for the equation
$\left( \ref{eqS14}\right) .$  Then $\left\|
V_1\left(t,s\right)\right\| =\left\| V\left( t,s\right) \right\|
_{1D}.$ Finally, applying  (\ref{eq2101}) we get the following
bound in $l_{1D}$ norm

\begin{eqnarray}
\sum\limits_{i\geq 1}q_i  p_i\left( t\right) =\left\| {\bf y}%
\left( t\right) \right\| \leq \left\| V_1\left( t,0\right) \right\|
\left\| {\bf y}\left( 0\right) \right\| +\int\limits_0^t\left\| V_1\left(
t,\tau \right) \right\|   \left\| D{\bf f}\left( \tau \right) \right\|
\,d\tau
\nonumber  \label{eqS15} \\
= \int\limits_0^t\left\| V_1\left(t,\tau \right) \right\|
  \left\| D{\bf f}\left( \tau \right) \right\|
\,d\tau \le
K\lambda_{0}e^{(-(l-\varepsilon) t} \int\limits_0^t
e^{\left((l-\varepsilon) \tau \right)} a(\tau) \, d\tau  \label{eqS15}
\\\le K\lambda_{0} a_{m} \left( 1+\frac{e^{(l-\varepsilon)T_{1}}}
{e^{(l-\varepsilon)T_{1}} - 1} \right).
\nonumber
\end{eqnarray}

Then
\begin{eqnarray*}
\sum\limits_{i>j}p_i\left( t\right) \leq q_{j+1} ^{-1}
\sum\limits_{i>j}q_i
p_i\left( t\right) \leq && \\
q_{j+1} ^{-1} K\lambda_{0} a_{m} \left( 1+\frac{e^{(l-\varepsilon)T_{1}}}
{e^{(l-\varepsilon)T_{1}} - 1} \right).
&&
\end{eqnarray*}
\vskip .9cm

Let  $E(t;k)=\sum\limits_{i\geq 1}i  p_i\left( t\right), \quad
t\ge 0, k\ge 0$ be the mean of the process at time $t,$ if
$X(0)=k.$

\begin{corollary}.
Under the assumptions of Theorem 1,

\begin{equation}
E( t;k) \leq \frac 1{{\rm W}} K\lambda_{0} a_{m} \left(
1+\frac{e^{(l-\varepsilon)T_{1}}} {e^{(l-\varepsilon)T_{1}} - 1}
\right), \forall \epsilon>0, \label{eqS17}
\end{equation}
where ${\rm W}=\inf\limits_{i\geq 1}\frac{q_i}i>0 $ and $l$ is
defined as in \refm[eq2101]
\end{corollary}

{\bf Proof.} The claim follows from the inequality
\[
\sum\limits_{i\geq 1}i  p_i\left( t\right) =\sum\limits_{i\geq 1}\frac
i{q_i}  q_i  p_i\left( t\right) \leq \frac 1{{\rm W}}
\sum\limits_{i \ge 1}q_i  p_i\left( t\right) ,\quad t\ge 0.
\]
\vskip .5cm

Consider now the case of finite state space.
In this case the matrix $D$ in (\ref{eqS7})
is also finite and we can obtain two-sided bounds
both for the decay function and for $\left\| {\bf p}^{(1)}( t) -
{\bf p}^{(2)}(t) \right\|,\quad t\ge 0.$

\begin{theorem}
Let the conditions of Theorem 1 hold and, in addition, $N <
\infty$.

Then

\nin (i) for any $s\geq 0,\ t\geq s,$ and any ${\bf %
p}^1\left( s\right) \in \Omega ,\ {\bf p}^2\left( s\right) \in \Omega $
\begin{eqnarray}
e^{-\int\limits_s^t\bar\zeta \left( u\right) du}
\left\| {\bf p}^1\left( s\right) -{\bf p}^2\left( s\right) \right\|_{1D} \le
\left\| {\bf p}^{(1)}\left( t\right) -{\bf p}^{(2)}\left( t\right) \right\|_{1D} \nonumber  \\
\le e^{-\int\limits_s^t l \left( u\right) du}
\left\| {\bf p}^{(1)}\left( s\right) -{\bf p}^{(2)}\left( s\right) \right\|_{1D}
, \quad 0\leq s\leq t,   \label{eqS1101}
\end{eqnarray}

\nin and

\begin{eqnarray}
\ \frac g{4NG}  e^{-\int\limits_s^t\bar\zeta \left( u\right) du}
\left\| {\bf p}^1\left( s\right) -{\bf p}^2\left( s\right) \right\|
\leq \left\| {\bf p}^1\left( t\right) -{\bf p}^2\left( t\right)
\right\|  \nonumber  \label{eqf6} \\
\ \leq
\frac{4NG}g
e^{-\int\limits_s^t\underline\alpha \left( u\right) du}
\left\| {\bf p}^1\left( s\right) -{\bf p}^2\left( s\right) \right\|
\label{eqf6}
\end{eqnarray}
where $g=\min d_k,\ G=\max d_k,$ and
$\underline\alpha \left( t\right), \quad \bar \zeta(t), \quad t\ge 0$
are defined in $\left( \ref{eqf3}\right) $ and
$\left( \ref{eq4119}\right) $ respectively;

\nin (ii) for any $s\geq 0,\ t\geq s,$  and
${\bf p}^1( s) \leq {\bf p}^2(s)$

\begin{equation}
\left\| {\bf p}^1\left( t\right) -{\bf p}^2\left( t\right) \right\|_{1D}
\geq e^{-\int\limits_s^t\bar \alpha \left( u\right)
du}
\left\| {\bf p}^1\left( s\right) -{\bf p}^2\left( s\right) \right\|_{1D} ,
\label{eqf1101}
\end{equation}

\nin and
\begin{equation}
\left\| {\bf p}^1\left( t\right) -{\bf p}^2\left( t\right) \right\|
\geq \frac g{4NG}  e^{-\int\limits_s^t\bar \alpha \left( u\right)
du}
\left\| {\bf p}^1\left( s\right) -{\bf p}^2\left( s\right) \right\| .
\label{eqf11}
\end{equation}
\end{theorem}

{\bf Proof.} We have to prove only estimates in $l_{1}$ norm. We
have

\begin{equation}
\left\| D\right\| =\sum\limits_{i=0}^{N-1}d_i\leq NG;\quad \left\|
D^{-1}\right\| =2  \max d_k^{-1}\leq \frac 2g.  \label{eqf7}
\end{equation}
Then

\begin{eqnarray}
&&\left\| {\bf p}^1\left( t\right) -{\bf p}^2\left( t\right) \right\|
\leq 2\left\| {\bf z}^1\left( t\right) -{\bf z}^2\left( t\right)
\right\| \leq  \nonumber  \label{eqf8} \\
&&2\left\| D\right\|   \left\| D^{-1}\right\|
e^{-\int\limits_s^t\underline \alpha \left( u\right) du}  \left\| {\bf
z}%
^1\left( s\right) -{\bf z}^2\left( s\right) \right\|  \label{eqf8} \\
&\leq &\frac{4NG}g  e^{-\int\limits_s^t\underline\alpha \left( u\right)
du}
\left\| {\bf z}^1\left( s\right) -{\bf z}^2\left( s\right) \right\|
\leq  \nonumber \\
&&\frac{4NG}g  e^{-\int\limits_s^t\underline\alpha
\left( u\right) du}
\left\| {\bf p}^1\left( s\right) -{\bf p}^2\left( s\right) \right\| .
\nonumber
\end{eqnarray}

On the other hand inequality $\left( \ref{eq1118}\right) $ implies
the bound
\begin{eqnarray}
&&\left\| {\bf p}^1\left( t\right) -{\bf p}^2\left( t\right) \right\|
\geq \left\| {\bf z}^1\left( t\right) -{\bf z}^2\left( t\right)
\right\| \geq  \nonumber  \label{eqf10} \\
&&\left\| D\right\| ^{-1}  \left\| D^{-1}\right\| ^{-1}
e^{-\int\limits_s^t\bar\zeta \left( u\right) du}  \left\| {\bf
z}^1\left(
s\right) -{\bf z}^2\left( s\right) \right\|   \\
&\geq &\frac g{2NG}  e^{-\int\limits_s^t\bar\zeta \left( u\right)
du}
\left\| {\bf z}^1\left( s\right) -{\bf z}^2\left( s\right) \right\|\geq  \nonumber \\
&&\frac g{4NG}  e^{-\int\limits_s^t\bar\zeta \left( u\right) du}
\left\| {\bf p}^1\left( s\right) -{\bf p}^2\left( s\right) \right\| .
\nonumber
\end{eqnarray}

\vskip .5cm

\begin{remark}. The rate of convergence in the previous theorems is exponential.
 Therefore, our bounds provide an  approach for
  obtaining integral estimates for the state probabilities and for the mean.
   Such estimates were obtained in \cite{doo00}, \cite{stad1}.
\end{remark}

\section{Null-ergodic BDPs}

\nin Put

\begin{equation}
\alpha_k^{(0)}\left( t\right)
=\lambda _k\left( t\right) +\mu _k\left( t\right)
-\delta _{k+1}  \lambda _k\left( t\right) -\delta _k^{-1}  \mu
_k\left( t\right) ,  \quad k=0,\ldots ,N-1, \quad t\ge 0.
\label{eq201,5}
\end{equation}
\nin and

\begin{equation}
\underline \alpha^{(0)} \left( t\right) =
\inf\limits_k\alpha  _k^{(0)}\left( t\right).
\label{eq202,5}
\end{equation}

Let ${\rm D}$ be a diagonal matrix
\begin{equation}
{\rm D}=diag\left\{ d_0,d_{1,\ldots }\right\},  \label{eqN1}
\end{equation}
where   $d_0=1,\quad d_i>0,\quad i\in E. $

\nin We define  the induced $l_{1{\rm D}}-$ norm  by

\[
 \left\| {\bf x}\right\| _{1{\rm %
D}}=\left\| {\rm D}{\bf x}\right\| _1=\sum\limits_{i\geq 0}d_i\left|
x_i\right|, \quad {\bf x}=(x_1,\ldots,\ldots )\in l_1.
\]
In the case considered all results are based on the bounds for
 the logarithmic norm $\gamma \left( A\left( t\right) \right) _{1{\rm D}%
}=\gamma \left( {\rm D}A\left( t\right) {\rm D}^{-1}\right) _1,$
where $A(t), \quad t \ge 0$ is defined in ( \ref{eq112}).

\nin Hence,
\begin{equation}
{\rm D}A(t){\rm D}^{-1}=\left(
\begin{array}{ccccc}
-\lambda _0(t) & d_0  d_1^{-1}  \mu _1(t) & 0 &  & \ldots \\
d_1  d_0^{-1}  \lambda _0(t) & -\left( \lambda _1(t)+\mu
_1(t)\right) &
d_1  d_2^{-1}  \mu _2(t) & 0 &  \\
0 & d_2  d_1^{-1}  \lambda _1(t) & \ddots & \ddots & \ddots \\
\vdots & 0 & \ddots & \ddots & \ddots \\
& \ldots & \ddots & \ddots & \ddots
\end{array}
\right), \quad t\ge 0.  \label{eqN2}
\end{equation}

Put
\begin{equation}
h_k=
\frac{\Delta \mu \lambda _{k-1}-\lambda \mu _k}{\Delta \lambda \mu
},\quad k\in \tilde E .
\label{eqn201}
\end{equation}

\begin{theorem}
Let for some $\Delta>1$
\begin{equation}
\inf_{k\in \tilde E}h_k = h > 0,  \label{eqn202}
\end{equation}
\nin let $c$ be a positive number under the conditions:
\begin{equation}
c < 1; \ c < \frac{\lambda_{k}}{\lambda}, k \ge 0, \ c \le h,  \label{eqn203}
\end{equation}
\nin and let
\begin{equation}
\lambda a_{m} - \Delta \mu b_{m}  > 0.
\label{eqn204}
\end{equation}
Then there exists a sequence $\mbox{\boldmath
$\delta$}=(\delta_k,\quad k\in \tilde E)>{\bf 0}$, such that

\nin (i)
\begin{equation}
\underline\alpha ^{(0)}\left( t\right)
\ge \theta(t) = c  \left( \lambda a\left( t\right)
-\Delta\mu b\left( t\right) \right),   \label{eqn205}
\end{equation}
\nin  where $\lim \inf \delta _k < 1, \quad {if} \quad N=\infty; $

\nin (ii) BDP $X(t),\quad t\ge 0$ is null-ergodic and
\begin{equation}
\sum\limits_{i=0}^\infty d_ip_i\left( t\right) \leq G
e^{-\int\limits_s^t\theta \left( u\right) du}, \quad 0\le s\le t,
\label{eqn206}
\end{equation}
where $G=\sup d_k<\infty .$

\end{theorem}

{\bf Proof.}

First we prove the inequalities

\begin{equation}
\alpha _k^{(0)}\left( t\right) =
a\left( t\right)   \lambda _k  \left(1-\delta _{k+1}\right)
-b\left( t\right)   \mu _k  \left( \delta
_k^{-1}-1\right) \geq c  \left( \lambda a\left( t\right) -\Delta\mu
b\left( t\right) \right),
 \label{eqn211}
\end{equation}

\nin $ k\in \tilde E,\quad t\ge 0. $

\nin It is sufficient to verify that

\begin{equation}\lambda _k  \left( 1-\delta _{k+1}\right) \geq c\lambda ,
\ \ \mu _k
\left( \delta _k^{-1}-1\right) \leq c\mu\Delta,\quad
 k\in \tilde E.
  \label{eqn212}
\end{equation}

\nin and (\ref{eqn212}) follows from

\begin{equation}
\delta _k\leq \frac{\lambda_{k-1}- c \lambda }{\lambda_{k-1}},\label{eqn213}
\end{equation}

\begin{equation}
\delta_k \ge \frac{\mu_k}{\mu_{k+1}+c\mu\Delta }, \quad k\in \tilde E.
\label{eqn214}
\end{equation}

If $N=\infty$ and $0<c<1,$ we have from (\ref{eqn213})

\begin{equation}
\limsup \delta_k \le \frac{\lambda -c\lambda}{\lambda} = 1-c < 1.
\label{eqn215}
\end{equation}

\nin Hence, (\ref{eqn213}), (\ref{eqn214}) imply

\begin{equation}
\frac{\mu_k}{\mu_{k+1}+c\mu\Delta } \le
\frac{\lambda_{k-1}- c \lambda }{\lambda_{k-1}}
, \quad k\in \tilde E  \label{eqn216}
\end{equation}

\nin Since $0<c< \frac{\lambda_{k}}{\lambda}$, the latter is
equivalent to
\begin{equation}
c^2\Delta\lambda \mu -c  \left( \Delta\lambda_{k-1} \mu -\lambda
\mu_{k} \right) \leq 0,
\quad k\in \tilde E.
\label{eqn216}
\end{equation}

\nin For fixed $k$ the largest root of $\left( \ref{eqn216}\right)
$ is $h_k$ given by \refm[eqn201].

\nin So, we have

\begin{equation}
\gamma \left( A\left( t\right) \right) _{1{\rm D}}\leq -\theta(t)
 = -c  \left( \lambda a\left( t\right) -\Delta\mu b\left(
t\right) \right),\quad t\ge 0
\label{eqn217}
\end{equation}

\nin Then
\begin{equation}
\sum\limits_{i=0}^\infty d_ip_i\left( t\right) =\left\| {\bf p}\left(
t\right) \right\| _{1{\rm D}}\leq e^{-\int\limits_s^t\theta
 \left(
u\right) du}  \left\| {\bf p}\left( s\right) \right\| _{1{\rm D}%
}\leq G  e^{-\int\limits_s^t\theta \left( u\right) du},
\quad 0\le s \le t,
\label{eqn218}
\end{equation}

\nin for any ${\bf p}\left( s\right) \in \Omega .$ This implies
the  exponential null-ergodicity of BDP.

\vskip .5cm

Now we are in position  to bound the state probabilities in
null-ergodic case. Put $d_j^{\min }=\min\limits_{i\leq j}d_i>0,
\quad j\in E.$

\begin{corollary}.
Under assumptions of Theorem 4,
\begin{equation}
p_k\left( t\right) \leq \frac G{d_k}  e^{-\int\limits_s^t\theta \left(
u\right) du}, \quad k\in E, \quad 0\le s\le t.
\label{eqn221}
\end{equation}
\nin and
\begin{equation}
\Pr \left( X\left( t\right) \leq j\mid X\left( 0\right) =k\right) \leq
\frac{%
d_k}{d_j^{\min }}  e^{-\int\limits_0^t\theta \left( u\right) du},
\quad j,k\in E,\quad t\ge 0
\label{eqn222}
\end{equation}

\end{corollary}

{\bf Proof.} We have
\begin{equation}
p_k\left( t\right) \leq d_k^{-1}  \sum\limits_{i=0}^\infty d_ip_i\left(
t\right) \leq \frac G{d_k}  e^{-\int\limits_s^t\theta \left( u\right)
du}.
\label{eqn223}
\end{equation}

\nin Thus,
\begin{eqnarray}
d_j^{\min }  \sum\limits_{i=0}^jp_i\left( t\right) \leq
\sum\limits_{i=0}^jd_ip_i\left( t\right) \leq \sum\limits_{i=0}^\infty
d_ip_i\left( t\right) &&  \nonumber \\
\leq e^{-\int\limits_0^t\theta \left( u\right) du}  \left\| {\bf p}%
\left( 0\right) \right\| _{1{\rm D}}=d_k  e^{-\int\limits_0^t\theta
\left( u\right) du}, &&  \label{eqn224}
\end{eqnarray}
which implies (\ref{eqn222}).
\vskip .5cm

\begin{remark} If BDP is null-ergodic then ${\rm E}\left( t;k\right)
\rightarrow \infty $ as $t\rightarrow \infty $ for any fixed $k\in E.$
Put
\begin{equation}
r\left( t\right) =\min \left\{ \lambda_0\left( t\right);\   \inf_{i\ge 1}
\left(
\lambda _i\left( t\right) -\mu _i\left( t\right) \right) \right \},\quad
t\ge 0.
\label{eqN8}
\end{equation}

Then $\frac{d{\rm E}\left( t;k\right) }{dt} \ge r(t)$, and

\begin{equation}
{\rm E}\left( t;k\right) \geq k+\int\limits_0^t r\left( \tau \right)
\,d\tau .  \label{eqN9}
\end{equation}

\end{remark}

\vskip .5cm

\section{ $M_t/M_t/1$ queue}

There is a  number of papers devoted to this queue, see the
references in \cite {mas85}. Here the number of customers $X\left(
t\right) $ at time $t$ (the length of the queue)  a BDP given by
\begin{equation}
\lambda _n\left( t\right) =a\left( t\right) ;\ \mu _n\left( t\right)
=b\left( t\right)  \label{eqMM1.1}
\end{equation}
 Bounds for this model were obtained in
\cite{z91a,z91b,z95a,z95b}.

Let $\rho =\frac{a_{m} }{b_{m} }$ be the traffic intensity.
 We consider underloaded ($\rho <1$) and overloaded ($\rho >1$) cases.

Put
\begin{equation}
\dot \alpha \left( t\right) =\left( 1-\rho ^{1/2}\right)   \left(
b\left( t\right) -\rho ^{-1/2}a\left( t\right) \right) .  \label{eqMM1.41}
\end{equation}

Put $\delta _i=\rho ^{-1/2},\ $ $d_i=\rho ^{-i/2},\ i\geq 1.$
Then we have in Theorem 1 $c=1-\rho ^{1/2},$ $g=\inf
d_k=1,$\ $g_i=\sum\limits_{m=0}^{i-1}d_m=\sum\limits_{m=0}^{i-1}\rho
^{-m/2}, \quad l ( t) =\dot \alpha \left( t\right) ,$ and
Theorem 1 implies the following bound

\begin{proposition}.
For underloaded $M_t/M_t/1$ queue ($a_{m} < b_{m}$) the following
upper bound on the rate of convergence holds:
\begin{equation}
\left\| {\bf p}^1\left( t\right) -{\bf p}^2\left( t\right) \right\|
\leq 4e^{-\int\limits_s^t\dot \alpha \left( u\right) du}
\sum\limits_{i\geq 1}\left\{ \sum\limits_{m=0}^{i-1}\rho ^{-m/2}\right\}
  \left| p_i^1\left( s\right) -p_i^2\left( s\right) \right| .
\label{eqMM1.5}
\end{equation}
\end{proposition}

\begin{remark}.
It is easy to see from \refm[eqMM1.41] that in the case considered
$\dot \alpha (t) =\left( \sqrt{a(t)}-\sqrt{b(t)}\right) ^2$ which
is the "exact" decay function, since the exact value of the
spectral gap (decay parameter) for the corresponding homogeneous
case is $\dot \alpha =\left( \sqrt{a}-\sqrt{b}\right) ^2$ , see
\cite{doorn85}.
\end{remark}

\begin{corollary}.
Let intensities be periodic with  period $T_{1} =T_{2} = 1$ and
$a_m=\int\limits_0^1 a(u) du < b_m =
\int\limits_0^1 b(u)du.$
Then  $\dot \alpha (t) $ is also $1$-periodic and
\begin{equation}
\dot \alpha_m:=\int\limits_0^1\dot \alpha(u) du=\left(
\sqrt{%
a_m}-\sqrt{b_m}\right)^2.  \label{eqMM1.6}
\end{equation}
\end{corollary}

\vskip .7cm

Theorem 2 and its Corollary imply the following statements

\begin{proposition}.
Let $M_t/M_t/1$ queue be underloaded. Then for any $\varepsilon >
0$ there exists $K,$ such that
\begin{equation}
\Pr \left( X\left( t\right) \leq j\mid X\left( 0\right) =0\right) \geq
\ 1-Ka_{m} \left\{ \sum\limits_{m=0}^j\rho ^{-m/2}\right\} ^{-1}
\left( 1+\frac{e^{(\dot \alpha -\varepsilon)T_{1}}}
{e^{(\dot \alpha -\varepsilon)T_{1}} - 1} \right).
\label{eqMM1.7}
\end{equation}
\end{proposition}

\begin{corollary}.
Let $M_t/M_t/1$ queue be underloaded. Then for any $\varepsilon >
0$ there exists $K,$ such that
\begin{equation}
{\rm E}\left( t;0\right) \le
\frac { Ka_{m}}{{\rm W}} \left( 1+\frac{e^{(\dot \alpha -\varepsilon)T_{1}}}
{e^{(\dot \alpha -\varepsilon)T_{1}} - 1} \right).
\label{eqMM1.8}
\end{equation}

where
\[
{\rm W}=\inf\limits_{i\geq 1}\left\{ i^{-1}
\sum\limits_{m=0}^{i-1}\rho ^{-m/2}\right\} >0.
\]
\end{corollary}

Consider the overloaded queue.

The conditions of Theorem 4 hold for $\delta _i=\rho ^{-1/2}<1,\ $
$d_i=\rho ^{-i/2},\ c=\rho ^{1/2} - 1.$ Then $G=1,\ \theta \left(
t\right) =\dot \alpha \left( t\right) $ and we obtain the
following statement

\begin{proposition}.
For overloaded $M_t/M_t/1$ queue, ($a_{m} > b_{m}$) the following
bound of the rate of convergence holds:
\begin{equation}
\sum\limits_{i=0}^\infty \rho ^{-i/2}p_i\left( t\right) \leq
e^{-\int\limits_s^t\dot \alpha \left( u\right) du},  \label{eqMM1.13}
\end{equation}
for any $s\geq 0,\ t\geq s.$

\end{proposition}

\begin{corollary}.
If $M_t/M_t/1$ queue is overloaded, then for any ${\bf p}\left(
s\right) \in \Omega $ and any $k$
\begin{equation}
p_k\left( t\right) \leq \rho ^{k/2}  e^{-\int\limits_s^t\dot \alpha
\left( u\right) du}.  \label{eqMM1.14}
\end{equation}
\end{corollary}

The following bound follows from Corollary of Theorem 4.

\begin{proposition}.
If $M_t/M_t/1$ queue is overloaded, then
\begin{equation}
\Pr \left( X\left( t\right) \leq j\mid X\left( 0\right) =k\right) \leq \rho
^{\left( k-j\right) /2}  e^{-\int\limits_0^t\dot \alpha \left( u\right)
du},  \label{eqMM1.15}
\end{equation}
\nin for any $t \ge 0$.
\end{proposition}

 With the help of  (\ref{eqN9}) we get the following bound for
 the mean of the  queue length.

\begin{proposition}.
\label{MM1.6} For any $k$ and any $t\geq 0$
\begin{equation}
{\rm E}\left( t;k\right) \geq k+\int\limits_0^t\left( a\left( \tau
\right) -b\left( \tau \right) \right) \,d\tau .  \label{eqMM1.16}
\end{equation}
\end{proposition}

\section{ $M_t/M_t/S$ queue}

For references see \cite{green91}, \cite{roth79}, \cite{stad1},
\cite{art}. Here the length of the queue (the number of customers)
is BDP $X\left(t\right) $ with the intensities
\begin{equation}
\lambda _n\left( t\right) =a\left( t\right) ;\ \mu _n\left( t\right)
=b\left( t\right)   \min \left( n,S\right).  \label{eqMMS.1}
\end{equation}

 Bounds for this process were obtained in \cite{z91a,z91b},
\cite{z95a}, \cite{z95b}. The estimates of probabilities for the
perturbed $M_t/M_t/S$ queue were studied in \cite{z02}.

Let $\rho =\frac{a_{m} }{S b_{m} }$ be the traffic intensity.
We consider underloaded ($\rho <1$) and overloaded ($\rho >1$) cases.

\begin{proposition}.
Let $M_t/M_t/S$ queue be underloaded. Then ( \ref{eq206}) holds
for
\begin{equation}
c=\min \left( \frac 1S,\frac{\Delta -1}\Delta \right), \quad \Delta>1 .  \label{eqMMS.101}
\end{equation}
\end{proposition}

\vskip .5cm

\begin{proposition}.
Let $M_t/M_t/S$ queue be overloaded. Then (\ref{eqn205}) holds for
\begin{equation}
c=\frac{\Delta -1}\Delta , \quad \Delta>1.  \label{eqMMS.102}
\end{equation}

\end{proposition}

\vskip .5cm

First, let the queue be underloaded. Put $\Delta = \rho ^{-1/2}$.
Then $c=\min \left( \frac 1S, 1 - \rho ^{1/2} \right)$. We
consider the cases of heavy traffic $\rho \geq
\left(\frac{S-1}S\right) ^2, \quad c= 1 - \rho ^{1/2},$ and light
traffic $0 \le \rho < \left(\frac{S-1}S\right) ^2, \quad c= \frac
1S .$

Put
\begin{equation}
\dot \alpha \left( t\right) = c   \left(
S b\left( t\right) -\rho ^{-1/2}a\left( t\right) \right) .  \label{eqMMS.2}
\end{equation}

Put $\delta _i=\rho ^{-1/2},\ $ $d_i=\rho ^{-i/2},\ i\geq 1.$ Then
we have in Theorem 1, $g=\inf d_k = 1$,
$q_i=\sum\limits_{m=0}^{i-1}d_m=\sum\limits_{m=0}^{i-1}
\rho^{-m/2}, \quad l( t) =\dot \alpha \left( t\right) ,$ and by
Theorem 1 we have

\begin{proposition}.
For underloaded $M_t/M_t/S$ queue ($a_{m} < S b_{m}$), the
following bound for the rate of convergence is valid:
\begin{equation}
\left\| {\bf p}^1\left( t\right) -{\bf p}^2\left( t\right) \right\|
\leq 4e^{-\int\limits_s^t\dot \alpha \left( u\right) du}
\sum\limits_{i\geq 1}\left\{ \sum\limits_{m=0}^{i-1}\rho ^{-m/2}\right\}
  \left| p_i^1\left( s\right) -p_i^2\left( s\right) \right| .
\label{eqMMS.3}
\end{equation}
\end{proposition}

\begin{remark}.
For underloaded $M_t/M_t/S$ queue with heavy traffic
  the exact "decay function" is
 $\dot \alpha (t) =\left( \sqrt{a(t)}-\sqrt{S b(t)}\right) ^2$
  because the exact value of spectral gap
  (decay parameter) for homogeneous case with heavy
  traffic is $\dot \alpha =\left( \sqrt{a}-\sqrt{Sb}\right) ^2$ ,
   see \cite{doorn85}.
     Moreover, for $S=2$ the bound between light and heavy traffics is also exact.
\end{remark}

\begin{corollary}.
Let intensities be periodic with period $T_{1} =T_{2} = 1$ and $a_m=\int\limits_0^1 a\left( u\right) du<S b_m=S \int\limits_0^1 b\left( u\right)
du.$ Then  $\dot \alpha (t) $ is also $1$-periodic and

\nin in the case of heavy traffic
\begin{equation}
\dot \alpha_m:=\int\limits_0^1\dot \alpha(u) du=
\left(\sqrt{a_m}-\sqrt{S b_m}\right)^2.
\label{eqMMS.6}
\end{equation}

\nin in the case of light traffic
\begin{equation}
\dot \alpha_m:=\int\limits_0^1\dot \alpha(u) du=
b_m -\sqrt{S a_m b_m}.
\label{eqMMS.601}
\end{equation}

\end{corollary}

\vskip .7cm

Theorem 2 and its Corollary imply the following two statements

\begin{proposition}.
Let $M_t/M_t/S$ queue be underloaded. Then for any $\varepsilon >
0$ there exists $K,$ such that
\begin{equation}
\Pr \left( X\left( t\right) \leq j\mid X\left( 0\right) =0\right) \geq
\ 1-Ka_{m} \left\{ \sum\limits_{m=0}^j\rho ^{-m/2}\right\} ^{-1}
\left( 1+\frac{e^{(\dot \alpha -\varepsilon)T_{1}}}
{e^{(\dot \alpha -\varepsilon)T_{1}} - 1} \right).
\label{eqMMS.7}
\end{equation}
\end{proposition}

\begin{corollary}.
Let $M_t/M_t/S$ queue be underloaded. Then for any $\varepsilon >
0$ there exists $K,$ such that
\begin{equation}
{\rm E}\left( t;0\right) \le
\frac { Ka_{m}}{{\rm W}} \left( 1+\frac{e^{(\dot \alpha -\varepsilon)T_{1}}}
{e^{(\dot \alpha -\varepsilon)T_{1}} - 1} \right).
\label{eqMMS.8}
\end{equation}

where
\[
{\rm W}=\inf\limits_{i\geq 1}\left\{ i^{-1}
\sum\limits_{m=0}^{i-1}\rho ^{-m/2}\right\} >0.
\]
\end{corollary}

Let now $M_t/M_t/S$ queue be overloaded. Put $\Delta = \rho
^{1/2}$. Then $c=  \rho ^{1/2} - 1 $. The conditions of Theorem 4
hold for $\delta _i=\rho ^{-1/2}<1,\ $ $d_i=\rho ^{-i/2},\
c=1-\rho ^{1/2}.$ Then $G=1,\ \theta \left( t\right) =\dot \alpha
\left( t\right) $ and we obtain

\begin{proposition}.
For overloaded $M_t/M_t/S$ queue ($a_{m} >S b_{m}$) the following
bound for the rate of convergence is valid:
\begin{equation}
\sum\limits_{i=0}^\infty \rho ^{-i/2}p_i\left( t\right) \leq
e^{-\int\limits_s^t\dot \alpha \left( u\right) du},
\label{eqMMS.13}
\end{equation}
for any $s\geq 0,\ t\geq s.$

\end{proposition}

\begin{corollary}.
If $M_t/M_t/S$ queue is overloaded, then for any ${\bf p}\left(
s\right) \in \Omega $ and any $k$
\begin{equation}
p_k\left( t\right) \leq \rho ^{k/2}  e^{-\int\limits_s^t\dot \alpha
\left( u\right) du}.  \label{eqMMS.14}
\end{equation}
\end{corollary}

The following bound follows from Corollary of Theorem 4.

\begin{proposition}.
If $M_t/M_t/S$ queue is overloaded, then
\begin{equation}
\Pr \left( X\left( t\right) \leq j\mid X\left( 0\right) =k\right) \leq \rho
^{\left( k-j\right) /2}  e^{-\int\limits_0^t\dot \alpha \left( u\right)
du},  \label{eqMMS.15}
\end{equation}
\nin for any $t \ge 0$.
\end{proposition}

>From (\ref{eqN9}) one gets the following simple bound for the mean
of the length of the queue

\begin{proposition}.
For any $k$ and any $t\geq 0$
\begin{equation}
{\rm E}\left( t;k\right) \geq k+\int\limits_0^t\left( a\left( \tau
\right) -S b\left( \tau \right) \right) \,d\tau .  \label{eqMMS.16}
\end{equation}
\end{proposition}

\vskip 2cm

\section{A multiserver queue with discouragements}

This queueing model was studied in
\cite{doorn81,nat74,re68,sri82}. We have here BDP with the
intensities $\lambda _n\left( t\right) =\lambda_n  a\left(
t\right) ;\ \mu _n\left( t\right) =\mu _n  b\left(t\right), $
where
\begin{equation}
\ \lambda _n=\left\{
\begin{array}{cc}
\lambda & \mbox {if }n<S \\
\frac \lambda {n-S+2} & \mbox {if }n\geq S
\end{array}
\right. ,\ \mu _n=\left\{
\begin{array}{cc}
n\mu & \mbox {if }n\leq S \\ S\mu & \mbox {if }n>S
\end{array}
\right. .  \label{eqdis.1}
\end{equation}

So, in the case considered $\lim\limits_{n\rightarrow \infty
}\lambda _n=0,$ which says that  a straightforward application of
Theorem 1 is impossible. Put
\begin{equation}
\delta _k=\left\{
\begin{array}{cc}
1 & \mbox {if }k<S \\
1+\varepsilon & \mbox {if }k\geq S
\end{array}
\right.  \label{eqdis.2}
\end{equation}
where $\varepsilon \in \left( 0;1\right) $. Then we have

\[
\alpha _k\left( t\right) =b\left( t\right), \quad  k<S-1;
\]
\[
\alpha _{S-1}\left( t\right) =b\left( t\right) -a\left( t\right)
\left( \frac{1+\varepsilon }2-1\right) ,
\]

\begin{equation}
\alpha _{n+S}\left( t\right) =S  b\left( t\right)   \left( 1-\frac
1{1+\varepsilon }\right) -a\left( t\right)   \left( \frac{1+\varepsilon
}{n+3}%
-\frac 1{n+2}\right), \quad  n=k-S \geq 0. \label{eqdis2.5}
\end{equation}
Put
\begin{equation}
c_n=\frac{\left( n+2\right) \varepsilon -1}{\left( n+2\right) \left(
n+3\right)
}  \label{eqdis3}
\end{equation}
If $\left( n+2\right) \leq \varepsilon ^{-1},$ then $c_n\leq 0;$ if $\left( n+2\right) >\varepsilon ^{-1},$ then $c_n<\frac\epsilon {n+3}<\frac{\varepsilon ^{2}}{1+\varepsilon} <\frac{S\varepsilon ^{2}}{1+\varepsilon}.$ Finally we obtain  for any $k$ the following bound
\begin{equation}
\alpha _k\left( t\right) \geq \dot \alpha (t) = \frac{S\varepsilon }{1+\varepsilon}   \left(
b\left( t\right) -\varepsilon   a\left( t\right) \right).  \label{eqdis4}
\end{equation}
Let $b_m > 0$ and let $\varepsilon \in \left( 0;1\right), $ be
such that $b_m - \varepsilon a_m > 0$. Then applying the estimates
of Theorem 1 we obtain

\begin{proposition}.
\begin{equation}
\left\| {\bf p}^1\left( t\right) -{\bf p}^2\left( t\right) \right\|
\leq 4  e^{-\int\limits_s^t\dot \alpha \left( u\right) du}
\sum\limits_{i\geq 1}g_i\left| p_i^1\left( s\right) -p_i^2\left( s\right)
\right| , \label{eqdis6}
\end{equation}
\nin where \ $q_i=\sum\limits_{m=0}^{i-1}\prod\limits_{k=0}^m\delta _k.$
\end{proposition}

Also,  Corollary 2 and Theorem 2 imply, under the same assumption
as above,

\begin{proposition}

Let $b_m > 0$. Then there exists $K,$ such that
\begin{equation}
\Pr \left( X\left( t\right) \leq j\mid X\left( 0\right) =0\right) \geq
1 - q_{j+1} ^{-1} K\lambda a_{m} \left( 1+\frac{e^{(\dot \alpha -\varepsilon)T_{1}}}
{e^{(\dot \alpha - \varepsilon)T_{1}} - 1} \right).
\label{eqdis7}
\end{equation}
\end{proposition}

\begin{corollary}.
Let $b_m > 0$. Then there exists $K,$ such that
\begin{equation}
{\rm E}\left( t;0\right) \le
\frac { Ka_{m}}{{\rm W}} \left( 1+\frac{e^{(\dot \alpha -\varepsilon)T_{1}}}
{e^{(\dot \alpha -\varepsilon)T_{1}} - 1} \right).
\label{eqdis8}
\end{equation}
\nin where ${\rm W}=\inf\limits_{i\geq 1}\left( q_i / i \right) >0. $

\end{corollary}

\section{ $M_t/M_t/S/S$ queue}

There is a large number of investigations of the $M/M /S/S$ queue,
see references in \cite{ki90}. In particular, this queue is a
subject of some new papers, see \cite{FRT}, \cite{voit},
\cite{stad2}. In general (nonstationary case) the length of a
queue $X\left(t\right) $ is a BDP on the state space $E=\left\{
0,1,\ldots ,S\right\} $ with the  intensities
\begin{equation}
\lambda _n\left( t\right) =a\left( t\right) ;\ \mu _n\left( t\right)
=nb\left( t\right)  \label{eqloss1}
\end{equation}
The first bounds for the general case were obtained in \cite{z89},
see also \cite{z95a,z95b}; some of results of \cite{z89} were
repeated in \cite{mas94}, and in the last paper non-Markovian case
was also studied.

Consider the condition
\begin{equation}
a_m + b_m > 0.  \label{eqloss2}
\end{equation}

One can see that BDP is not exponentially ergodic if ( \ref{eqloss2}) is not fulfilled.

Let ( \ref{eqloss2}) holds.

{\bf The first case.}
\nin Let
\begin{equation}
b_m > 0.  \label{eqloss3}
\end{equation}

Put $\delta _k=1,\ k=1,\ldots ,S-1.$ Then
$ d_k = 1,\ g=\min d_k = 1,\ G=\max d_k = 1,$ and
\begin{equation}
\alpha _k\left( t\right) =\lambda _k\left( t\right) +\mu _{k+1}\left(
t\right) -\delta _{k+1}  \lambda _{k+1}\left( t\right) -\delta
_k^{-1}  \mu _k\left( t\right)  =
\left\{
\begin{array}{cc}
b\left( t\right) & {if } k<S-1 \\
b\left( t\right) +a\left( t\right) & {if } k=S-1
\end{array}
\right.
\label{eqloss4}
\end{equation}
\nin and
\begin{equation}
\tilde \alpha \left( t\right) =\inf\limits_k\alpha _k\left( t\right)
=b\left( t\right) ,  \label{eqloss5}
\end{equation}

\begin{equation}
\tilde \beta \left( t\right) =\sup\limits_k\alpha _k\left( t\right)
=a\left(
t\right) +b\left( t\right),  \label{eqloss6}
\end{equation}

\begin{equation}
\tilde \zeta \left( t\right) =\sup\limits_k\left\{ \lambda _k\left(
t\right) +\mu _{k+1}\left( t\right) +\delta _{k+1}  \lambda
_{k+1}\left(
t\right) +\delta _k^{-1}  \mu _k\left( t\right) \right\}
\leq 2a\left( t\right) +\left( 2S-1\right) b\left( t\right) .
\label{eqloss7}
\end{equation}

Now Theorem 3 implies

\begin{proposition}.
Let the condition ( \ref{eqloss3}) be satisfied. Then for any
 $s\geq 0,\ t\geq s,$ ${\bf p}^1\left( s\right) \in \Omega
,\ {\bf p}^2\left( s\right) \in \Omega $

\begin{equation}
e^{-\int\limits_s^t\tilde \zeta \left( u\right) du}
\left\| {\bf p}^1\left( s\right) -{\bf p}^2\left( s\right) \right\|_{1D} \le
\left\| {\bf p}^{(1)}\left( t\right) -{\bf p}^{(2)}\left( t\right) \right\|_{1D}
\le e^{-\int\limits_s^t \tilde \alpha \left( u\right) du}
\left\| {\bf p}^{(1)}\left( s\right) -{\bf p}^{(2)}\left( s\right) \right\|_{1D}
,   \label{eqloss701}
\end{equation}

\nin and

\begin{equation}
\frac 1{4S}  e^{-\int\limits_s^t\tilde \zeta \left( u\right) du}
\left\| {\bf p}^1\left( s\right) -{\bf p}^2\left( s\right) \right\|_1
\leq \left\| {\bf p}^1\left( t\right) -{\bf p}^2\left( t\right)
\right\|_1
 \leq 4S  e^{-\int\limits_s^t\tilde \alpha \left( u\right) du}
\left\| {\bf p}^1\left( s\right) -{\bf p}^2\left( s\right) \right\|_1 .
\label{eqloss8}
\end{equation}
\end{proposition}
\vskip.7cm

\begin{remark}.
If the intensities are constant, then the bounds (\ref{eqloss701})
and (\ref{eqloss8}) immediately imply the inequality (17) in the
paper \cite{FRT}. It is interesting to note that such bound in
general case was found about ten years ago in \cite{z89}. The
cutoff property, which is studied in \cite{FRT} and \cite{voit}
can be naturally investigated by our methods. It will be a subject
of a separate paper.
\end{remark}

\begin{proposition}.
Let the condition ( \ref{eqloss3}) be satisfied. Then for any
 $s\geq 0,\ t\geq s,$ and ${\bf p}^1\left( s\right) \leq {\bf p}^2\left( s\right) $ the following bounds hold

\begin{equation}
\left\| {\bf p}^1\left( t\right) -{\bf p}^2\left( t\right) \right\|_{1D}
\geq e^{-\int\limits_s^t\tilde \beta \left( u\right)
du}
\left\| {\bf p}^1\left( s\right) -{\bf p}^2\left( s\right) \right\|_{1D} ,
\label{eqloss801}
\end{equation}

\nin and

\begin{equation}
\left\| {\bf p}^1\left( t\right) -{\bf p}^2\left( t\right) \right\|_1
\geq \frac 1{4S}  e^{-\int\limits_s^t\tilde \beta \left( u\right)
du}  \left\| {\bf p}^1\left( s\right) -{\bf p}^2\left( s\right)
\right\|_1 .  \label{eqloss9}
\end{equation}
\end{proposition}
\vskip.7cm

Consider now the bounds for the length of queue. We have here ${\rm W}=\inf\limits_{i\geq 1}\frac{g_i}%
i=\inf\limits_{i\geq 1}\frac 1i  \left\{
\sum\limits_{m=0}^{i-1}\prod\limits_{k=0}^m\delta _k\right\} =1.$

Theorem 2 and its Corollary imply the following statetment.

\begin{proposition}
Let the condition ( \ref{eqloss3}) be satisfied. Then for any
$\varepsilon > 0$ there exists $K,$ such that
\begin{equation}
\Pr \left( X\left( t\right) \leq j\mid X\left( 0\right) =0\right) \geq
1 - \frac{ K a_{m}}{j+1}  \left( 1+\frac{e^{(b_m-\varepsilon)T_{1}}}
{e^{(b_m-\varepsilon)T_{1}} - 1} \right).
\label{eqloss91}
\end{equation}

\nin and

\begin{equation}
{\rm E}\left( t;0\right) \leq
K a_{m} \left( 1+\frac{e^{(b_m-\varepsilon)T_{1}}}
{e^{(b_m-\varepsilon)T_{1}} - 1} \right).
\label{eqloss92}
\end{equation}

\end{proposition}

On the other hand, the estimate of the mean of the number of
customers in the queue can be obtained in the following  way:
\begin{eqnarray}
\ &&\frac{d{\rm E}\left( t;k\right) }{dt}=\lambda _0\left( t\right)
p_0+\sum\limits_{i=1}^S\left( \lambda _i\left( t\right) -\mu _i\left(
t\right) \right)   p_i=a\left( t\right)
\sum\limits_{i=0}^{S-1}p_i-b\left( t\right)   \sum\limits_{i=0}^Si
p_i\leq  \nonumber  \label{eqN911} \\
\ &&a\left( t\right)   \sum\limits_{i=0}^Sp_i-b\left( t\right)
{\rm E}\left( t;k\right) =a\left( t\right) -b\left( t\right)
{\rm E}\left( t;k\right) .  \label{eqloss93}
\end{eqnarray}
 By integration of this inequality we get

\begin{proposition}
The following bound for the mean of the length of  $M_t/M_t/S/S$ queue holds:
\begin{equation}
{\rm E}\left( t;k\right) \leq k  e^{-\int\limits_0^tb\left( u\right)
du}+\int\limits_0^ta\left( \tau \right)   e^{-\int\limits_\tau
^tb\left(
u\right) du}d\tau  \label{eqloss94}.
\end{equation}
\end{proposition}

\begin{remark}.
The bound (\ref{eqloss94}) seems to be sufficiently sharp for the case of light traffic. It is interesting to study its sharpness in the general case.
\end{remark}
\vskip.7cm

{\bf The second case.}
Let now
\begin{equation}
a_m > 0.
\label{eqloss10}
\end{equation}

Put $\delta _k=\delta =\frac{S-1}S<1,\ k=1,\ldots ,S-1.$ Then
\begin{eqnarray}
\alpha _k\left( t\right) &=&\lambda _k\left( t\right) +\mu _{k+1}\left(
t\right) -\delta _{k+1}  \lambda _{k+1}\left( t\right) -\delta
_k^{-1}  \mu _k\left( t\right)  \label{eqloss11} \\
&=&\left\{
\begin{array}{cc}
\left( 1-\delta \right)   a\left( t\right) +\left( k-\frac{k-1}\delta
\right) b\left( t\right) & {, if  }\  k<S-1 \\
a\left( t\right) & {, if  } \  k=S-1
\end{array}
\right.  \nonumber
\end{eqnarray}
\nin and
\begin{equation}
\breve \alpha \left( t\right) =\inf\limits_k\alpha _k\left( t\right) \geq
\left( 1-\delta \right)   a\left( t\right) =\frac{a\left( t\right) }S,
\label{eqloss12}
\end{equation}

\begin{equation}
\breve \beta \left( t\right) =\sup\limits_k\alpha _k\left( t\right) \leq
a\left( t\right) +b\left( t\right) =\tilde \beta \left( t\right),
\label{eqloss13}
\end{equation}

\begin{equation}
\breve \zeta \left( t\right) =\sup\limits_k\left\{ \lambda _k\left(
t\right)
+\mu _{k+1}\left( t\right) +\delta _{k+1}  \lambda _{k+1}\left(
t\right)
+\delta _k^{-1}  \mu _k\left( t\right) \right\} \leq \left( 1+\delta
\right)   a\left( t\right) +2S  b\left( t\right) .
\label{eqloss14}
\end{equation}

Then by virtue of Theorem 3, we have

\begin{proposition}
Let the condition ( \ref{eqloss10}) hold. Then for any $s\geq 0,\
t\geq s,$ and any ${\bf p}^1\left( s\right) \in \Omega,\ {\bf
p}^2\left( s\right) \in \Omega $

\begin{equation}
e^{-\int\limits_s^t \breve \zeta \left( u\right) du}
\left\| {\bf p}^1\left( s\right) -{\bf p}^2\left( s\right) \right\|_{1D} \le
\left\| {\bf p}^{(1)}\left( t\right) -{\bf p}^{(2)}\left( t\right) \right\|_{1D}
\le e^{-\int\limits_s^t \breve \alpha \left( u\right) du}
\left\| {\bf p}^{(1)}\left( s\right) -{\bf p}^{(2)}\left( s\right) \right\|_{1D}
,   \label{eqloss1401}
\end{equation}

\nin and

\begin{eqnarray}
&&\frac{S-1}{4S^2}  e^{-\int\limits_s^t\breve \zeta \left( u\right)
du}  \left\| {\bf p}^1\left( s\right) -{\bf p}^2\left( s\right)
\right\|_1 \leq \left\| {\bf p}^1\left( t\right) -{\bf p}^2\left(
t\right) \right\|_1  \label{eqloss15} \\
&\leq &\frac{4S^2}{S-1}  e^{-\int\limits_s^t\breve \alpha \left(
u\right) du}  \left\| {\bf p}^1\left( s\right) -{\bf p}^2\left(
s\right) \right\|_1 .  \nonumber
\end{eqnarray}

\end{proposition}
\vskip0.5cm

\begin{proposition}
Let the condition ( \ref{eqloss10}) hold. Then for any
 $s\geq 0,\ t\geq s,$ and ${\bf p}^1\left( s\right) \leq {\bf p}^2\left( s\right) $
 the following bounds are valid

\begin{equation}
\left\| {\bf p}^1\left( t\right) -{\bf p}^2\left( t\right) \right\|_{1D}
\geq e^{-\int\limits_s^t \breve \beta \left( u\right)
du}
\left\| {\bf p}^1\left( s\right) -{\bf p}^2\left( s\right) \right\|_{1D} ,
\label{eqloss1501}
\end{equation}

\nin and

\begin{equation}
\left\| {\bf p}^1\left( t\right) -{\bf p}^2\left( t\right) \right\|_1
\geq \frac{S-1}{4S^2}  e^{-\int\limits_s^t\breve \beta \left( u\right)
du}  \left\| {\bf p}^1\left( s\right) -{\bf p}^2\left( s\right)
\right\|_1 .  \label{eqloss16}
\end{equation}
\end{proposition}

\nin {\bf Acknowledgement}

 The research of the first  author was supported by the Fund for
 the Promotion of Research at Technion.

\end{document}